{\rm }\input amstex
\documentstyle{amsppt}

\def\ms{\medskip}
\nologo
\def\un{\underbar}

\def\nl{\newline}
\def\l[{\leftbracket}
\def\r]{\rightbracket}

\NoBlackBoxes
\NoRunningHeads
\magnification=\magstep1
\pageheight{6.8in}
\def\tstrut{\vrule height 3.1ex depth 1.2ex width 0pt}
\TagsOnRight
\topmatter
\title A DOUBLE BOUNDED VERSION OF SCHUR'S PARTITION THEOREM\endtitle
\author { Krishnaswami Alladi and Alexander Berkovich}\endauthor
\affil University of Florida, Gainesville, Florida  32611\\
\ms
{\it{Dedicated to the memory of Professor Paul Erd\"os}} \endaffil
\subjclass Primary -05A15, 05A17, 11PP81, 11P83\endsubjclass
\keywords partitions, Schur's theorem, method of weighted words, key 
identity, Chu-Vandemonde summation, q-binomial coefficients\endkeywords
\abstract
Schur's partition theorem states that the number of partitions of $n$
into distinct parts $\equiv 1, 2$ (mod $3$) equals the number of partitions
of $n$ into parts which differ by $\geqslant 3,$ where the inequality is
strict if a part is a multiple of 3.  We establish a double bounded
refined version of this theorem by imposing one bound on the parts 
$\equiv 0,1$ (mod $3$) and another on the parts $\equiv 2$(mod $3$), and by 
keeping track of the number of parts in each of the residue classes
(mod 3). Despite the long history of Schur's theorem, our result is
new, and extends earlier work of Andrews, Alladi-Gordon and
Bressoud.  We give combinatorial and q-theoretic proofs of our
result. The special case L=M leads to a representation of the generating 
function of the underlying partitions in terms of the q-trinomial 
coefficients extending a similar previous representation of Andrews.
\endabstract
\endtopmatter
\document
\baselineskip=17pt
\head \S1. {\un{Introduction}}\endhead
Schur's celebrated partition theorem of 1926 is the following result
[11]: 
\proclaim{Theorem S} The number of partitions of $n$ into distinct
parts $\equiv1, 2$ (mod $3$) equals the number of partitions of $n$ into
parts differing by $\geqslant 3,$ where consecutive multiples of $3$
cannot occur as parts. 
\endproclaim

Subsequently Gleissberg [9] extended Theorem S to an arbitrary modulus
$m\geqslant 3,$ and established a stronger correspondence involving the
number of parts.

Several proofs of Schur's theorem by a variety of approaches are
known, most notably by Andrews [4] using generating functions, by
Bressoud [7] involving a combinatorial bijection, and by Andrews
[5] using a representation involving the q-trinomial coefficients. 
Alladi and Gordon [3] obtained generalizations and refinements using 
a new technique called {\it{the method of weighted words}}. They viewed a
strong refinement of Schur's theorem as emerging out of the {\it{key
identity}}

$$\sum_{i,j,\ge 0} A^iB^j \sum_{k=0}^{\min(i,j)} \dfrac
{q^{T_{i+j-k}+T_k}}{(q)_{i-k} (q)_{j-k} (q)_{k}} =
\sum_{i,j,\geqslant 0}
\dfrac{A^iB^iq^{Ti+Tj}}{(q)_i(q)_j}=(-Aq)_\infty (-Bq)_\infty\tag1.1$$
under the transformations
$$\matrix\text {(dilation)} q\longmapsto q^3,\\
\text{(translations)} A\longmapsto Aq^{-2},\  B\longmapsto Bq^{-1}.
\endmatrix\bigg\}\tag1.2$$
In (1.1) and in what follows, $T_n=n(n+1)/2$ is the n-th triangular
number, and the symbols $(A)_n$ are defined by
$$(A)_n=(A;q)_n=
\cases
\prod^{n-1}_{j=0}(1-Aq^j), \qquad if \qquad n>0,\\
1, \qquad \qquad \qquad \qquad if \qquad n=0,\\
\prod^{-n}_{j=1}{(1-Aq^{-j})^{-1}}, \, if \qquad n<0.
\endcases
\tag1.3$$
when $n$ is an integer, and 
$$(A)_\infty={\lim_{n\to\infty}}(A)_n=\prod_{j=0}^{\infty}(1-Aq^j), 
\qquad {\text{for}} \quad |q|<1.\tag1.4$$ 

Our goal here is to establish the finite identity (2.1) stated in
\S2, from which (1.1) follows when L, M$\rightarrow\infty$.  When the
transformations (1.2) are applied to (2.1), the combinatorial
interpretation yields Theorem 3 in \S 4, which is the double bounded
refined version of Theorem S, with different bounds on the parts
$\equiv 0,1$ (mod $3$) and $\equiv 2$ (mod $3$), and where the number of
parts in the three residue classes can be specified.  Although there is a 
long and rich history on Schur's theorem, and refinements keeping track of
number of parts in residue classes are known (such as what follows
from the infinite identity (1.1)), our result with bounds on
parts in residue classes never seems have been stated in the
literature.

In \S2 we give a proof of the finite identity (2.1) by making use of the 
q-Chu-Vandermonde summation. Then we give two more proofs of (2.1) with 
the conditions $L,M\ge i+j$ and $i,j\ge 0$, in which case the identity has 
combinatorial interpretation. More precisely, the second proof in \S3
utilizes Durfee rectangles and extends the ideas in [3], and the third
proof which is combinatorial and bijective (see \S4) uses the methods 
in [3] and [7]. 

The case $L=M$ in (2.1) is discussed in \S5. Alladi and Gordon [3] had shown 
that the refined Schur theorem underlying (1.1) emerged from the study 
of the numerator of a certain infinite continued fraction. In \S5 we show 
that our generating function corresponding to the bound $L$ on the parts 
are the numerator convergents to this continued fraction. In addition this 
also leads to a representation involving the q-trinomial coefficients 
(see \S5) with two free parameters $A$ and $B$ in the summation, extending 
a result of Andrews [5] who had previously obtained such a representation 
with one free parameter.

The method of weighted words was substantially improved by Alladi, Andrews, 
and Gordon [1], to obtain generalizations and refinements of a deep partition 
theorem of G\"ollnitz. They did this by proving the infinite key identity 
(6.1) which is considerably deeper than (1.1). We have recently obtained a 
double bounded finite version of this G\"ollnitz key identity. This is stated 
as identity (6.3) in \S6 and a complete discussion of it will be taken up 
later [2].

\head \S2. {\un{A double bounded key identity}}\endhead
Let $L, M, i, j,$ be arbitrary integers. Then we have 
$$\sum_{k=0}^{min(i, j)} q^{(i-k)(j-k)}\left[\matrix{M-i-j+k}\\k\endmatrix\right]\left[\matrix{M-j}\\i-k\endmatrix\right]\left[\matrix{L-i}\\j-k\endmatrix\right]=\left[\matrix{L}\\j\endmatrix\right]\left[\matrix{M-j}\\i\endmatrix\right],\tag2.1$$
\nl where the $q$-binomial coefficients 
$\left[\matrix{n+m}\\{n}\endmatrix\right]$, are defined by
$$
\left[\matrix{n+m}\\{n}\endmatrix\right]_q =
\left[\matrix{n+m}\\{n}\endmatrix\right] =
\cases
\frac{(q^{m+1})_n}{(q)_n}, \qquad if \quad n\ge0,\\
0, \qquad \qquad \quad if \quad n<0.
\endcases\tag2.2
$$ 
As in (2.2), when the base for the binomial coefficient is $q$, we will  
supress it, but if the base is anything other than $q$, it will be displayed. 
When $m\ge0$, (2.2) yields the standard definition of the q-binomial 
coefficient which is symmetric in m and n, namely, 
$$
\left[\matrix{n+m}\\{m}\endmatrix\right] = 
\cases
\dfrac{(q)_{n+m}}{(q)_m(q)}_n, \qquad if \quad n\ge0,\\
0, \qquad \qquad \quad if \quad n<0.
\endcases\tag2.3
$$ 

Let us now show that (2.1) yields (1.1) when $L,M\rightarrow\infty$. From 
(2.2) it follows that as $L, M\rightarrow\infty$, identity (2.1) reduces to
$$
\sum_{k=0}^{\min(i,j)}\dfrac{q^{(i-k)(j-k)}}{(q)_k(q)_{i-k}(q)_{j-k}}
=\dfrac{1}{(q)_i(q)_j}.\tag2.4
$$
Next it is easy to verify that 
$$
T_{i+j-k} + T_k = T_i + T_j + (i-k)(j-k).\tag2.5
$$
Thus (2.4) and (2.5) yield
$$
\sum_{k=0}^{\min(i,j)}\dfrac {q^{T_{i+j-k}+T_k}}{(q)_{i-k}(q)_{j-k}(q)_k}=
\dfrac{q^{T_i+T_j}}{(q)_i(q)_j},\tag2.6
$$
which is (1.1), by comparing the coefficients of $A^iB^j$ on both
sides.

Now we give a proof of (2.1).  To this end we use (I.10) and (I.20) in Gasper 
and Rahman [8] to rewrite the left hand side of (2.1) as 
$$
\lim_{l\to L} \frac{(q^{l-i-j+1})_j(q^{M-i-j+1})_i}{(q)_j(q)_i}q^{ij}
\sum_{k\ge0}\frac{(q^{-j})_k(q^{-i})_kq^k}{(q)_k(q^{l-i-j-1})_k},\tag2.7
$$
where we used limit definitions to make sure that all objects in (2.7) are well
defined. The sum in (2.7) can be evaluated by appeal to the
q-Chu-Vandermonde summation formula (II.6 in [8]) as 
$$
\frac{(q^{l-j+1})_jq^{-ij}}{(q^{l-i-j+1})_j}.\tag2.8
$$
Combining (2.7) and (2.8) we see that the left hand side of (2.1) is 
$$
\lim_{l\to L} \frac{(q^{M-i-j+1})_i(q^{l-j+1})_j}{(q)_i(q)_j}
= \left[\matrix{M-j}\\i\endmatrix\right]\left[\matrix{L}\\j\endmatrix\right]
$$
which is the right hand side of (2.1) completing the proof.

\head \S3.{\un{Proof using Durfee rectangles.}}\endhead

We now give another proof of (2.1) when $L,M\ge{i+j}$ and $i,j\ge 0$. In 
this case, we can use (2.3) to write (2.1) in the fully expanded form 
$$
\sum_{k=0}^{\min(i,j)} q^{(i-k)(j-k)}\dfrac{(q)_{M-i-j+k}}{(q)_k(q)_{M-i-j}}
\qquad\dfrac{(q)_{M-j}}{(q)_{i-k}(q)_{M-i-j+k}}\qquad
\dfrac{(q)_{L-i}}{(q)_{j-k}(q)_{L-i-j+k}}
$$
$$= \dfrac {(q)_L}{(q)_j(q)_{L-j}} \cdot 
\dfrac {(q)_{M-j}}{(q)_i(q)_{M-j-i}}.\tag3.1
$$
Note that all terms involving the parameter $M$ disappear from (3.1)
after cancellation!  Thus in this case, (2.1) is equivalent to the identity
$$
\sum_{k=0}^{\min(i,j)}q^{(i-k)(j-k)}\left[\matrix i\\k
\endmatrix\right]\left[\matrix L-i\\j-k\endmatrix\right] =
\left[\matrix L\\j\endmatrix\right]\tag3.2
$$
which can be proved combinatorially using Durfee rectangles.

The term $\left[\matrix L\\j\endmatrix\right]$ on the right hand side of
(3.2) is the generating function of partitions $\pi$ into $\leqslant
j$ parts with each part $\leqslant L-j$.  In the Ferrers graph of each
such partition $\pi$, there is a maximal Durfee rectangle whose row length
minus column length is $i-j$.  Let this rectangle have $j-k$ rows and
$i-k$ columns.  The number of nodes in this rectangle is $(i-k)(j-k)$
and this accounts for the term 
$$
q^{(i-k)(j-k)}
$$
in (3.2).  Next, the portion of the Ferrers graph below this rectangle
is a partition into no more than $k$  parts each $\leqslant i-k$.  The
generating function of such partitions is 
$$
\left[\matrix i\\k\endmatrix\right]
$$
Finally the partition of the Ferrers graph to the right of this
rectangle is a partition into $\leqslant j-k$ parts each $\leqslant
L-i-j+k.$  The generating function of such partitions is
$$
\left[\matrix L-i\\j-k\endmatrix\right]
$$
Thus
$$
q^{(i-k)(j-k)} \left[\matrix i\\k\endmatrix\right]\left[\matrix
L-i\\j-k\endmatrix\right]
$$
is the generating function of all such Ferrers graphs $\pi$ having a
Durfee rectangle of size $(i-k)(j-k)$. The parameter $k$ depends on
the particular Ferrers graph given. We need to sum over $k$ to
account for all Ferrers graphs under discussion, and so this proves
(3.2).

\un{Remark:} This argument using Durfee rectangles is the same as
in Alladi-Gordon [3] except that they did not impose a bound $L-j$ on
the size of the parts.

\head \S4.{\un{Combinatorial proof.}}\endhead
We now give a combinatorial proof of (2.1) by following the method of
Alladi-Gordon [3]. In order to do this we need to briefly describe the
generalization of Schur's partition theorem established in [3].

Alladi and Gordon consider integers occurring in three colors, of
which two are primary represented by $a$ and $b,$ and one is secondary,
represented by $ab$.  The integer 1 occurs only in the primary colors $a$ and 
$b$, and the integers $n \geqslant 2$ occur in all three colors, $a, b,$
and $ab$.  By the symbols $ a_n, b_n, ab_n$, we mean the integer $n$
occurring in colors $a, b$, and $ab$, respectively.  To discuss partitions
into colored integers we need an ordering among the symbols and the
one we choose is 
$$
a_1<b_1<ab_2<a_2<b_2<ab_3<a_3<b_3...\tag4.1
$$
The reason for the choice of this ordering will be explained soon.

By a partition of $n$ we mean a sum of symbols arranged in decreasing
order according to (4.1) such that the sum of the subscripts (weights)
in $n$.  For example, $b_5+b_5+a_5+(ab)_5+a_4+b_3+(ab)_3$ is a
partition of 30. By a Type-1 partition we mean a partition of the form
$n_1+n_2+...+n_\nu,$ where the $n_i$ are symbols in (4.1) such that the
gap between the subscripts of consecutive symbols $n_i$ and $n_{i+1}$ in
this partition is $\geqslant 1$, but with strict inequality if
$$
\matrix n_i \text{ is of color } ab\text { or if } \\
n_i \text { is of color } a \text { and } n_{i+1} \text { is of color
} b\endmatrix\bigg\}.\tag4.2
$$
In [3] it is shown that 
$$
\dfrac{q^{T_{i+j-k}+T_k}}{(q)_{i-k}(q)_{j-k}(q)_k}
$$
is the generating function of all Type-1 partitions involving
exactly $i-k$ $a$-parts, $j-k$ $b$-parts, and $k$ $ab$-parts. The term
$$
\dfrac{q^{T_i+T_j}}{(q)_i(q)_j}
$$
is clearly the generating function for all vector partitions $(\pi_1,
\pi_2)$ of $N$ where $\pi_1$ has $i$ distinct a-parts and $\pi_2$ has
$j$ distinct b-parts.  In view of these explanations, the combinatorial 
version of (1.1) is the following result [3].
\proclaim{Theorem 1}
Let V(n;i,j) denote the number of vector partitions $(\pi_1,\pi_2)$ of $N$ 
where $\pi_1$ has exactly $i$ distinct a-parts and $\pi_2$ has exactly $j$ 
distinct b-parts.  

Let $S(n;r,s,t)$ denote the number of Type-1 partitions of $n$ having 
$r$ $a$-parts, $s$ $b$-parts, and $t$ $ab$-parts. 

Then
$$
V(n; i,j) = \sum\Sb r+t=i\\s+t=j\endSb S(n;r,s,t).
$$
\endproclaim

A strong refinement of Schur's theorem follows from Theorem 1 by using
the transformations (1.2) which converts the product on the right in
(1.1) to
$$
\overset\infty \to{\underset m=1\to\Pi} (1+Aq^{3m-2})(1+Bq^{3m-1}).
$$
Under the same transformations, the symbols $a_n, b_n, ab_n$ become 
$$
\matrix a_n\mapsto 3n-2,\quad b_n \mapsto 3n-1, \text { for
} n\geqslant 1,\\
ab_n \mapsto 3n-3, \text { for } n \geqslant 3.\endmatrix\bigg\}\tag4.3
$$
Thus under (4.3), the ordering (4.1) is just the natural ordering
among the positive integers
$$
1<2<3<...,
$$ 
which explains the reason for the choice of this ordering.  Also under
the transformations (4.3), the gap condition (4.2) governing the
Type-1 partitions become the Schur gap conditions, namely, the gap
between consecutive parts is $\geqslant3,$ with strict inequality if a
part is a multiple of 3.  Thus Theorem 1 is a strong refinement of
Schur's theorem in the undilated form.

In order to make the combinatorial proof of (2.1) as clear as possible
we need to give the combinatorial proof of Theorem 1 in [3] here.
Once that is done, we can go through the steps of that proof by
imposing bounds $L$ and $M$ on certain parts and then (2.1) will fall
out easily.

The combinatorial proof of Theorem 1 will be illustrated with the
vector partition $(\pi_1;\pi_2)$, where
$$
\pi_1: a_6+a_5+a_3+a_2+a_1,\qquad \pi_2:b_9+b_8+b_6+b_4+b_2+b_1
$$
Here $i=5$ and $j=6.$

Suppose in general that $\pi_1$ has $i$ parts and $\pi_2$ has j parts.

{\un{Step 1:}} Decompose $\pi_2$ into $\pi_4$ and $\pi_5$, where $\pi_4$ 
has the parts of $\pi_2$ which are $\leqslant i$ and $\pi_5$ has the 
remaining parts:
$$
\pi_4: b_4+b_2+b_1\qquad\qquad \pi_5: b_9+b_8+b_6$$

{\un{Step2:}} Consider the conjugate of the Ferrers graph of
$\pi_4$ and circle the bottom node of each column.  Denote this graph
by $\pi_4^*.$  Construct a graph $\pi_6$ where the number of nodes in
each row is the sum of the number of nodes in the corresponding rows
of $\pi_1$ and $\pi_4^*.$  The parts of $\pi_6$ ending in circled nodes
are $ab$-parts.  The rest are $a$-parts:
$$
\pi_6 = \pi_1 + \pi_4^*: ab_9+ab_7+a_4+ab_3+a_1.
$$

Conversely, given $\pi_6$, the columns ending with the circled nodes 
can be extracted to form $\pi_4^*$, and what remains after the extraction 
will be $\pi_1$. 

{\un{Step 3:}} Write the parts of $\pi_5$ in a column in descending
order and below them write the parts of $\pi_6$ in descending order.

{\un{Step 4:}} Subtract $0$ from the bottom element, $1$ from the
next element above, $2$ from the one above that, etc., and display the
new values as well as the subtracted ones in two adjacent columns
$C_1|C_2.$ The elements of $C_2$ have no color, while those of $C_1$
retain the colors of the parts from which they were derived.
  
{\un{Step 5:}} Rearrange the elements of $C_1$ in decreasing order
given by (4.1) to form a column $C_1^R.$

{\un{Step 6:} Finally, add the corresponding elements of $C_1^R$ and
$C_2$ to get a partition $\pi_3$ counted by $S(n; r, s, t).$  The
colors of the parts of $\pi_3$ are those of the elements of $C_1^R$
from which they were derived.

\input tables
\tableinfofalse
\noncenteredtables
\line{
\begintable
Step 3\crthick
$\pi_5/\pi_6$\crthick
$b_9$\crnorule
$b_8$\crnorule
$b_6$\crnorule
$ab_9$\crnorule
$ab_7$\crnorule
$a_4$\crnorule
$ab_3$\crnorule
$a_1$
\endtable
\qquad\qquad\qquad
\begintable
\multispan{2}\tstrut\hfil Step 4\hfil\crthick
$C_1$ | $C_2$\crthick
$  b_2$| 7\crnorule
$ b_2$ | 6\crnorule
$ b_1$| 5\crnorule
$ab_5$ | 4\crnorule
$ab_4$| 3\crnorule
$ a_2$ | 2\crnorule
$ ab_2$ | 1\crnorule
$ a_1$ | 0
\endtable
\qquad\qquad\qquad
\begintable
\multispan{2}\tstrut\hfil Step 5\hfil\crthick
$C_1^R$ | $C_2^R$\crthick
$ab_5$| 7\crnorule
$ab_4$ | 6\crnorule
$ b_2$| 5\crnorule
$b_2$ | 4\crnorule
$a_2$| 3\crnorule
$ ab_2$ | 2\crnorule
$ b_1$ | 1\crnorule
$ a_1$ | 0
\endtable
\qquad\qquad\qquad
\begintable
Step 6\hfil\crthick
$\pi_3$\crthick
$ab_{12}$\crnorule
$ab_{10}$\crnorule
$b_7$\crnorule
$b_6$\crnorule
$a_5$\crnorule
$ab_4$\crnorule
$b_2$\crnorule
$a_1$\endtable
}

Each of these steps is a one-to-one correspondence, and so this is a
bijective proof of Theorem 1.

We now give a bijective proof of (2.1) when $L,M\ge {i+j}$ and $i,j\ge 0$ 
by discussing the above steps in reverse. For this purpose we need to 
consider the quantities $\nu(\pi;M)=\nu(M)$ and $\nu(\pi;L)=\nu(L)$, which 
which are defined as follows.
 
1) If $L \geq M$ then $\nu(L) = 0$. If $L < M$ then $\nu(L)$ is the
non-negative number such that there are exactly $\nu(L)$ parts in the interval
$[L-\nu(L)+2,M]$, all $b$ parts $\leq L-\nu(L)$ and no part $= L-\nu(L)+1$. 

2) If $M \geq L$ then $\nu(M) = 0$. If $M < L$ then $\nu(M)$ is the
non-negative number such that there are exactly $\nu(M)$ parts in the interval
$[M-\nu(M)+2,L]$, all $ab,a$ parts $\leq M-\nu(M)$ and no part $= M-\nu(M)+1$.

First using (2.5) we rewrite (2.1) in the equivalent form 
$$
\sum_{k=0}^{\min(i,j)}q^{T_{i+j-k}+T_k}\left[\matrix
M-i-j+k\\k\endmatrix\right] \left[\matrix
M-j\\i-k\endmatrix\right]\left[\matrix L-i\\j-k\endmatrix\right] =
q^{T_i+T_j}\left[\matrix L\\j\endmatrix\right]\left[\matrix
M-j\\i\endmatrix\right].\tag 4.4
$$
Identity (4.4) can be proved combinatorially using the above steps as
we show now.

First, we interpret
$$
q^{T_k}\left[\matrix M-i-j+k\\k\endmatrix\right]\left[\matrix
M-j\\i-k\endmatrix\right]\left[\matrix L-i\\j-k\endmatrix\right]
$$

as the generating function for partitions having 
$$
\matrix \leqslant i-k\quad  a\text {-parts, each }\leqslant M-j-i+k,\\
\leqslant j-k\quad b\text { -parts, each } \leqslant L-i-j+k,\\
\text { and } k \quad \text{distinct} \quad ab\text { -parts, each } 
\leqslant M-i-j+k.
\endmatrix\Biggr\}
$$
Next, to this partition, we add 1 to the smallest part, 2 to the second 
smallest part, ..., $i+j-k$ to the largest part, and retain the colors of 
the parts we started with. We then get
$$
q^{T_{i+j-k}{+T_k}}\left[\matrix
M-j-i+k\\k\endmatrix\right]\left[\matrix
M-j\\i-k\endmatrix\right]\left[\matrix L-i\\j-k\endmatrix\right]
$$
as the generating function for Type-1 partitions having
$$
\matrix i-k \quad a\text { -parts } \leqslant M - \nu(M),\\
j-k \quad b\text {-parts }\leqslant L-\nu(L),\\
k \quad ab\text { -parts } \leqslant M-\nu(M).\endmatrix\Biggr\}
$$
This is the interpretation of the summand on the left in (4.4) and the
partition corresponds to one of the form $\pi_3$ in Step 6.

Going from Step 6 to Step 5, we subtract 0, 1, 2, 3, ..., $i+j-k-1$ in
succession.  Thus the partitions in column $C_1^R$ would have:
$$
\matrix i-k \quad a\text {-parts each } \leqslant M-i-j+k+1,\\
j-k \quad b\text {-parts each } \leqslant L-i-j+k+1,\\
k \quad \text {distinct } ab \text {-parts each } \leqslant
M-i-j+k+1.\endmatrix\Biggr\} $$
Proceeding from Step 5 to Step 3, we do a rearrangement, and then we add
0, 1, 2, ..., $i-1$ to the $a$-parts and $ab$-parts, whereas we add $i,
i+1, ..., i+j-k-1$ to the $b$-parts.  Thus $\pi_5/\pi_6$ in Step 3 is a
partition with 
$$
\matrix i-k\quad \text {distinct } a\text { -parts each } \leqslant M-j+k,\\
j-k\quad \text {distinct } b\text { -parts in the interval } \left[i+1,
L\right],\\ k\quad ab\text { -parts differing by } \ge 2 \text { and all }
\leqslant M-j+k.\endmatrix\Biggr\} $$

Finally in Step 2 the $k$ $ab$-parts decompose into $k$ $a$-parts and
$k$ $b$-parts, the latter being $\leqslant i.$ In this process
all $a$-parts are bounded by $M-j.$ Thus we end up with partitions
having
$$
\matrix i \text { distinct } a\text {-parts each } \leqslant M-j,\\
\text { and } j \text { distinct } b\text { -parts each } \leqslant
L.\endmatrix\bigg\}\tag4.5
$$
The generating function of the partitions satisfying (4.5) is
$$
q^{T_i+T_j}\left[\matrix M-j\\i\endmatrix\right]\left[\matrix L\\j
\endmatrix\right]
$$ 
which is the right hand side of (4.4), thereby completing the
combinatorial proof.

In view of the above proof and interpretation, we can improve Theorem 1 
to the following double bounded form:
\proclaim{Theorem 2} 

Let $L, M, i, j$ be non-negative integers with $M\ge L\ge{i+j}$.

Let $V(n; i, j, L, M)$ denote the number of vector partitions
$(\pi_1;\pi_2)$ of $n$ having $i$ distinct a-parts each $\leqslant
M-j,$ and $j$ distinct b-parts each $\leqslant L.$

Let $S(n; r, s, t, l, L, M)$ denote the number of Type-1
partitions of $n$ having  $r$ $a$-parts $\leqslant M$, $s$ $b$-parts 
$\leqslant L-l$, $t$ $ab$-parts $\leqslant M$, there are exactly $l=\nu(L)$
$\text{ } a,ab$-parts which are $\geq L-l+2$, and no part $= L-l+1$.  

Then
$$ 
V(n;i,j,L,M) = \sum\Sb r+t=i\\ s+t=j \endSb\sum\Sb l\endSb S(n;r,s,t,l,L,M).
$$
\endproclaim

NOTE: Since $M\ge L$ in Theorem 2, we have $\nu(M)=0$ and so the inner 
summation is over $l=\nu(L)$. If $L\ge M$, then $\nu(L)=0$ and so there 
would be similar theorem with the inner summation over $m=\nu(M)$. 

Under the transformations (1.2), Theorem 2 yields the following double
bounded strong refinement of Schur's theorem which is new:

\proclaim{Theorem 3}  Let $L, M, i, j,$ be non-negative integers with 
$M\ge L\ge {i+j}$.  

Let $ P(n;i,j,L,M)$ denote the number of partitions of $n$ into $i$ distinct 
parts $\equiv 1$ (mod $3$) each $\leqslant 3(M-j)-2$, and $j$ distinct parts
$\equiv 2$ (mod $3$) each $\leqslant 3L-1.$

Let $G(n;r, s, t, l, L, M)$ denote the number of partitions of $n$ into
parts differing by $\geqslant 3,$ where the inequality is strict if a
part is a multiple of $3$, and such that there are $r$ parts $\equiv 1$ 
(mod $3$) each $\leqslant3M-2, s$ parts $\equiv 2$ (mod $3$) each
$\leqslant 3(L-l)-1$, t parts $\equiv 0$ (mod $3$) each $\leqslant 3M-3$, 
there are exactly $l$-parts which are $>3(L-l)+2$, and no part $= 3(L-l),
3(L-l)+1$. 

Then
$$
P(n; i, j, L, M) = \sum\Sb r+t=i\\s+t=j\endSb \sum\Sb l \endSb 
G(n; r, s, t, l, L,M).
$$
\endproclaim

As in the case of Theorem 2, there is a version of Theorem 3 when 
$L > M\ge i+j$. 

{\bf{Remarks}}: There are other ways in which one might obtain double 
bounded versions of the identity (2.6). For instance, consider 
$$
(-Aq)_M(-Bq)_L=\sum_{i=0}^{M}\sum_{j=0}^{N}A^iB^jq^{T_i+T_j}
\left[\matrix M\\i\endmatrix\right]\left[\matrix L\\j\endmatrix\right]. 
\tag4.6
$$
We may use the q-binomial theorem to expand $(-Aq)_M$ and rewrite the left 
hand side of (4.6) as
$$
(-Aq)_M(-Bq)_L=\sum_{i=0}^{M}A^iq^{T_i}(-Bq)_i(-Bq^{i+1})_{L-i}
\left[\matrix M\\i\endmatrix\right]. 
\tag4.7
$$
By expanding the factors $(-Bq)_i$ and $(-Bq^{i+1})_{L-i}$ using the 
q-binomial theorem once again, and rearranging, we get on comparison with 
(4.6) the following finite identity for Schur's theorem:
$$
q^{T_i+T_j}\left[\matrix M\\i\endmatrix\right]\left[\matrix L\\j
\endmatrix\right]=\sum \Sb k \endSb q^{T_{i+j-k}+T_k}
\left[\matrix M\\M-i,i-k,k\endmatrix\right]
\left[\matrix L-i\\j-k\endmatrix\right], \tag4.8
$$
where 
$$
\left[\matrix M\\i,j,M-i-j\endmatrix\right]=\frac{(q)_M}{(q)_i(q)_j(q)_{M-i-j}}
$$
is the q-multinomial coefficient of order 3. Identity (4.8) has the advantage 
that the bounds on the parts enumerated by the left hand side are simple, 
namely, that the $a$-parts are bounded by $M$ and the $b$-parts by $L$. But 
then the bounds on the parts of the partitions enumerated by 
the right hand are more complicated. The decomposition considered 
in (4.7) corresponds precisely to the decomposition of the partition $\pi_2$ 
into $\pi_4+\pi_5$ in Step 1 above.

\head \S5. {\un{Other Connections}}\endhead
In this section we discuss the case $L=M$ of the double bounded key
identity (4.4), namely, the case where all parts of $S(n; r, s, k)$ are
bounded by $b_L.$  In this case the product of the q-binomial
coefficients on the right hand side of (4.4) becomes
$$
\left[\matrix L\\j\endmatrix\right] \left[\matrix
L-j\\i\endmatrix\right] = \dfrac {(q)_L}{(q)_j(q)_{L-j}} \dfrac
{(q)_{L-j}}{(q)_i(q)_{L-i-j}} = \left[\matrix L\\i,j,L-i-j\endmatrix\right],
\tag5.1
$$
the q-multinomial coefficient (of order 3). Thus when $L=M$ (4.4) reduces to 
$$
\sum_{k=0}^{min(i,j)}q^{T_{i+j-k}+T_k}\left[\matrix
L-i-j+k\\k\endmatrix\right]\left[\matrix
L-j\\i-k\endmatrix\right]\left[\matrix L-i\\j-k\endmatrix\right] =
q^{T_i+T_j} \left[\matrix L\\ i,j, L-i-j\endmatrix\right].\tag5.2
$$

Now multiply both sides of (5.2) by $A^iB^j$ and sum over $i$ and $j$ to get
$$
\sum_{i,j\ge 0}A^iB^j\sum_{k=0}^{\min(i,j)}q^{T_{i+j-k}+T_k}
\left[\matrix
L-i-j+k\\k\endmatrix\right]\left[\matrix
L-j\\i-k\endmatrix\right]\left[\matrix L-i\\j-k\endmatrix\right]
$$
$$= \sum_{i,j\ge 0} A^iB^j q^{T_i+T_j}\left[\matrix
L\\i,j,L-i-j\endmatrix\right].\tag5.3
$$
Even though (5.3) is an immediate consequence of (4.4), it is instructive 
to give an independent proof of (5.3). We will do so by interpreting the 
left hand side of (5.3) combinatorially, and by utilizing a not so well 
known recurrence formula for the q-multinomial coefficients for the right 
hand side (see (5.8) below). This leads to a connection with a continued 
fraction expansion considered by Alladi and Gordon [3] and extends a 
representation due to Andrews [5] involving the q-trinomial coefficients. 

Denote the left hand side of (5.3) by $G_L (A, B; q)$. From the discussion 
in \S 4 it follows that $G_L(A,B;q)$ is the generating function of Type-1 
partitions $\pi$ where all parts are $\le b_L$. That is
$$
G_L (A, B; q) =\sum\Sb \pi \text { of type } 1\\ \lambda
(\pi)\leqslant b_L \endSb A^{\nu_a(\pi)}B^{\nu_b(\pi)}(AB)^{\nu_{ab}(\pi)}
q^{\sigma (\pi)},\tag5.4
$$
where $\sigma(\pi)$ is the sum of the parts of $\pi, \lambda (\pi)$ the
largest part of $\pi$, and $\nu_a(\pi)$, $\nu_b(\pi)$, $\nu_{ab}(\pi)$, denote
the number of $a$-parts, $b$-parts, and $ab$-parts of $\pi$ respectively.

In Alladi and Gordon [3] it is shown that the generating function 
$G_L(A,B;q)$ satisfies the recurrence 
$$
G_L(A, B; q) = (1+Aq^L + Bq^L) G_{L-1} (A, B; q) + ABq^L (1-q^{L-1})
G_{L-2} (A,B;q),\tag5.5
$$
but they did not have the representation for $G_L(A,B;q)$ as the left hand 
side of (5.3). 

Next, let $R_L (A, B; q)$ denote the right hand side of (5.3).  We
claim that $R_L$ satisfies the same recurrence, namely,
$$
R_L(A,B;q)=(1+Aq^L+Bq^L)R_{L-1}(A,B;q)+ABq^L(1-q^{L-1})R_{L-2}(A,B;q).\tag5.6
$$
By considering the coefficient of $A^iB^j$ in $R_L,$ we see that
proving (5.6) is equivalent to showing
$$ 
q^{T_i+T_j}\left[\matrix L \\ i,j,L-i-j\endmatrix\right] =
q^{T_i+T_j}\left[\matrix L-1 \\ i,j,L-i-j-1\endmatrix\right] + 
q^{T_{i-1}+T_j+L}\left[\matrix L-1 \\i-1, j, L-i-j\endmatrix\right]
$$
$$
+q^{T_i+T_{j-1}+L}\left[\matrix L-1
\\i,j-1,L-i-j\endmatrix\right]+q^{T_{i-1}+T_{j-1}+L}(1-q^{L-1})\left[\matrix
L-2 \\i-1, j-1, L-i-j\endmatrix\right].\tag5.7
$$
By cancelling $q^{T_i+T_j}$ on both sides of (5.7), we get the following
equivalent second order (in L) recurrence relation for the q-multinomial
coefficients:
$$
\left[\matrix L \\ i,j,L-i-j\endmatrix\right]=\left[\matrix
L-1\\i,j,L-1-i-j\endmatrix\right] +q^{L-i} \left[\matrix
L-1 \\ i-1,j,L-i-j\endmatrix\right]
$$
$$+ q^{L-j} \left[\matrix L-1\\i,j-1,L-i-j\endmatrix\right]
+q^{L-i-j}(1-q^{L-1})\left[\matrix L-2\\i-1,j-1,L-i-j\endmatrix\right].
\tag5.8
$$
This recurrence is symmetric in $i$ and $j$ but is not so well known;
it can be derived from one of the six (non-symmetric) standard
recurrences for the q-multinomial coefficients, namely,
$$
\left[\matrix L\\i,j,L-i-j\endmatrix\right]=\left[\matrix
L-1\\i,j,L-1-i-j\endmatrix\right]+q^{L-i-j}\left[\matrix
L-1\\i,j-1,L-i-j\endmatrix\right]
$$
$$
+q^{L-i}\left[\matrix L-1\\i-1,j,L-i-j\endmatrix\right]\tag5.9
$$
Indeed (5.9) implies that (5.8) is equivalent to
$$q^{L-i-j}\left[\matrix
L-1\\i,j-1,L-i-j\endmatrix\right]=q^{L-j}\left[\matrix
L-1\\i,j-1,L-i-j\endmatrix\right]
$$
$$
+q^{L-i-j}(1-q^{L-1})\left[\matrix L-2\\i-1,j-1,L-i-j\endmatrix\right]
\tag5.10
$$
We may cancel the common factor $q^{L-i-j}$ and $(q)_{L-1}$ in the
numerator of (5.10) and the common factors $(q)_{j-1}$ and
$(q)_{L-i-j}$ in the denominator and reduce (5.10) to 
$$\dfrac {1}{(q)_i} = \dfrac {q^i}{(q)_i} + \dfrac {1}{(q)_{i-1}}$$
which is clearly true.  Thus (5.8) is established, and consequently
(5.6).  Finally, since $G_L$ and $R_L$ satisfy the same initial
conditions, (5.3) is proven.

Alladi and Gordon [3] viewed the left hand side of (1.1) as the
numerator of
$$
\left(1+(A+B)q+\dfrac{ABq^2(1-q)}{1+(A+B)q^2+\dfrac{ABq^3(1-q^2)}{1+(A+B)q^3 
+ ...}}\right).\tag5.11
$$
Thus the infinite key identity is the statement that the numerator of
this continued fraction equals
$$\overset\infty \to{\underset m=1\to\Pi}(1+Aq^m)(1+Bq^m)$$
Now let $P_L(A, B; q)$ denote the numerator of the $L-th$ convergents
to this continued fraction.  Clearly $P_L$ satisfies the recurrence
$$
P_L(A, B; q) = (1 + Aq^L + Bq^L) P_{L-1} (A, B; q) + ABq^L
(1-q^{L-1}) P_{L-2} (A, B; q).\tag5.12$$  
By comparing (5.12) with (5.5) and checking initial conditions it
follows that the numerator convergents $P_L$ are precisely the
generating functions $G_L$ for Type-1 partitions $\pi$ with
$\lambda (\pi) \leq b_L.$ The left hand side of (5.3) is the
representation for the convergents $P_L$.  When $L\to\infty$, this
representation yields the left hand side of (1.1).

We conclude this section by producing a representation for $G_L$ in
terms of the q-trinomial coefficients, thereby extending a previous
similar representation by Andrews [5]. To this end, use the
transformations (1.2) to recast (5.3) in the form
$$
\sum_{i,j}\sum_{k=0}^{\min(i,j)}A^iB^jq^{3(T_{i+j-k}+T_k)-2i-j}\left[\matrix
L-i-j+k\\k\endmatrix\right]_{q3}\left[\matrix
L-i\\j-k\endmatrix\right]_{q3}\left[\matrix
L-j\\i-k\endmatrix\right]_{q^3}
$$
$$=\sum_{i,j}A^iB^jq^{3(T_i+T_j)-2i-j}\left[\matrix
L\\i,j,L-i-j\endmatrix\right]_{q^3}.\tag5.13
$$
Next, replace $i$ by $j+\tau$ on the right hand side of (5.13) to rewrite
it as 
$$\sum_{j,\tau}(Aq)^{j+\tau}(Bq^2)^j q^{3\left(\matrix
j+\tau\\2\endmatrix\right)+3\left(\matrix j\\2\endmatrix\right)}
\left[\matrix L\\j+\tau, j,L-2j-\tau\endmatrix\right]_{q^3}$$
$$=\sum_{j,\tau}(AB)^j
A^{\tau}q^{\tau\left(\dfrac{3\tau-1}{2}\right)+3j(j+\tau)}\left[\matrix L\\j+\tau,j,L-2j-\tau
\endmatrix\right]_{q^3}.\tag5.14$$
Following Andrews and Baxter [6] we define the generalized q-trinomial
coefficients$\left(\matrix L;q\\\tau\endmatrix\right)_c$ by 
$$
\left(\matrix L;q\\\tau\endmatrix\right)_c = \sum_j c^j
q^{j(j+\tau)}\left[\matrix L\\j+\tau,j,L-2j-\tau\endmatrix\right].\tag5.15
$$
Then from (5.3), (5.4), (5.13), (5.14) and (5.15), we get the representation
$$
G_{3L-1}(Aq^{-2},Bq^{-1};q^3)=\sum_\tau A^{\tau}q^{\tau\left(\dfrac{3\tau-1}{2}\right)}
\left(\matrix L;q^3 \\ \tau \endmatrix\right)_{c=AB}\tag5.16$$
for the generating function of all partitions enumerated by $B(n)$
having parts $\leq 3L-1.$ Previously Andrews [5] had utilized the 
$q$-trinomial coefficients to get the representation (5.16) where $c=AB=1$ 
which corresponds to Gleissberg's refinement [9] of Schur's partition 
theorem. 

\head \S6. {\un{A double bounded G\"ollnitz identity}}\endhead

One of the deepest results in the theory of partitions is a theorem of 
G\"ollnitz [10]:
\proclaim{Theorem G}  Let $P(n)$ denote the number of partitions of
$n$ into distinct parts $\equiv 2, 4$ or $5$ (mod $6$).

Let $G(n)$ denote the number of partitions of $n$ into parts $\neq 1$, 
or $3$, such that the difference between the parts is $\geqslant 6$ with 
strict inequality if a part is $\equiv 0, 1$ or $3$ (mod $6$).

Then
$$
G(n)=P(n).
$$
\endproclaim

The proof by G\"ollnitz [10] is very involved.

Alladi, Andrews, and Gordon [1], obtained substantial generalizations 
and refinements of Theorem G, and viewed this theorem as emerging out of the 
{\it{key identity}}
$$
\sum_{i,j,k}A^iB^jC^k\sum\Sb i=\alpha+\delta+\varepsilon\\ 
j=\beta+\delta+\phi\\ k=\gamma+\varepsilon+\phi\endSb
\frac{q^{T_s+T_{\delta}+T_{\varepsilon}+T_{\phi-1}}(1-q^{\alpha}+
q^{\alpha+\phi})} 
{(q)_{\alpha}(q)_{\beta}(q)_{\gamma}(q)_{\delta}(q)_{\varepsilon}(q)_\phi}
=(-Aq)_\infty(-Bq)_\infty(-Cq)_\infty\tag6.1
$$
under the transformations
$$
\matrix\text {(dilation)}  q\rightarrow q^6,\\
\text {(translations)} A\rightarrow Aq^{-4}, B\rightarrow Bq^{-2},
C\rightarrow Cq^{-1}.\endmatrix\bigg\}\tag6.2
$$
Here $s=\alpha+\beta+\delta+\delta+\varepsilon+\phi$. Note that the key 
identity (1.1) for Schur's theorem is the special case of (6.1) with $C=0$.

We have recently obtained the following double bounded version of (6.1): 

If $i,j,k,L,M$ are given integers, then  
$$
\align
&\sum q^{T_s + T_{\delta} + T_{\varepsilon} + T_{{\phi}-1}}\\
&\qquad \left\{q^{\phi}
\left[\matrix L-s+\alpha\\ \alpha\endmatrix\right] 
\left[\matrix L-s+\beta\\ \beta\endmatrix\right] 
\left[\matrix M-s+\gamma\\ \gamma\endmatrix\right] 
\left[\matrix L-s\\ \delta\endmatrix\right]
\left[\matrix M-s\\ \varepsilon\endmatrix\right]
\left[\matrix M-s\\ \phi\endmatrix\right]\right.\\
&\qquad \left. +  
\left[\matrix L-s+\alpha-1\\ \alpha-1\endmatrix\right] 
\left[\matrix L-s+\beta\\ \beta\endmatrix\right] 
\left[\matrix M-s+\gamma\\ \gamma\endmatrix\right] 
\left[\matrix L-s\\ \delta\endmatrix\right]
\left[\matrix M-s\\ \varepsilon\endmatrix\right]
\left[\matrix M-s\\ \phi-1\endmatrix\right] 
\right\} \\
\endalign
$$
$$
=\sum_{\tau\ge0} q^{\tau(M+2)-T_\tau+T_{i-\tau}+T_{j-\tau}+T_{k-\tau}}
\left[\matrix L-\tau\\\tau\endmatrix\right]
\left[\matrix L-2\tau\\i-\tau\endmatrix\right]
\left[\matrix L-i-\tau\\j-\tau\endmatrix\right]
\left[\matrix M-i-j\\k-\tau\endmatrix\right],
\tag6.3
$$
where $s$ is as in (6.1) and the summation in (6.3) is over $\alpha, \beta, 
\gamma, \delta, \varepsilon, \phi$ satisying the same conditions with 
respect to $i,j,k$ as in (6.1). In (6.3) if we set either $i=0$ or $j=0$, 
then we get an identity equivalent to (4.4). If we set $k=0$, then (6.3) 
reduces to (5.3) which is (4.4) with $L=M$. By letting the parameters 
$L,M\rightarrow\infty$, only the term corresponding to $\tau=0$ on the right 
in (6.3) survives, and  (6.3) reduces to 

$$
\sum\Sb i=\alpha+\delta+\varepsilon\\ 
        j=\beta+\delta+\phi\\
        k=\gamma+\varepsilon+\phi\endSb
\dfrac
{q^{T_s+T_{\delta}+T_{\varepsilon}+T_{\phi-1}}(1-q^{\alpha}+q^{\alpha+\phi})}
{(q)_{\alpha}(q)_{\beta}(q)_{\gamma}(q)_{\delta}(q)_{\varepsilon}(q)_{\phi}}
=\dfrac{q^{T_i+T_j+T_k}}{(q)_i(q)_j(q)_k}
$$
from which (6.1) follows if we multiply both sides by $A^iB^jC^k$ and sum 
over $i,j,k$. The proof of the new identity (6.3) will be the subject 
of a forthcoming paper [2].

Finally, if $L=M$, then the sum on the right in (6.3) can be evaluated 
using the q-Pfaff-Saalsch\"utz summation (see [8], formula II.12) in terms  
of three q-binomial coefficients where $i,j,k$ occur cyclically
$$
q^{T_i+T_j+T_k}\left[\matrix L-k\\i\endmatrix\right]
\left[\matrix L-i\\j\endmatrix\right]\left[\matrix L-j\\k\endmatrix\right].
$$
In this case (6.3) can be interpreted in terms of partitions (see [2]).

\enddocument